\documentclass[12pt]{article}
\usepackage{amssymb}
\usepackage{amsfonts}
\usepackage{latexsym}
\usepackage[dvips]{graphics}

\newtheorem{thm}{Theorem}
\newtheorem{lem}{Lemma}

\newtheorem{cor}{Corollary}
\newtheorem{prop}{Proposition}
\newtheorem{rem}{Remark}

\newtheorem{ex}{Example}

\large\normalsize

\begin{document}
\begin{center}
{\Large\bf Replacing Pfaffians and applications}
\\[30pt]
{Weigen\ Yan$^{\rm a,b}$ \footnote{This work is supported by
FMSTF(2004J024) and NSFF(E0540007)} \quad and \quad Yeong-Nan\
Yeh$^{\rm b}$ \footnote{Partially supported by NSC94-2115-M001-017
\newline\hspace*{5mm}{\it Email address:} weigenyan@263.net (W.
Yan), mayeh@math.sinica.edu.tw (Y. N. Yeh)}}
\\[10pt]
\footnotesize { $^{\rm a}$School of Sciences, Jimei University,
Xiamen 361021, China
\\[7pt]
$^{\rm b}$Institute of Mathematics, Academia Sinica, Taipei 11529,
Taiwan}
\end{center}
\begin{abstract}
We present some Pfaffian identities, which are completely
different from the Pl\"ucker relations. As consequences we obtain
a quadratic identity for the number of perfect matchings of plane
graphs, which has a simpler form than the formula by Yan et al
(Graphical condensation of plane graphs: a combinatorial approach,
Theoret. Comput. Sci., to appear), and we also obtain some new
determinant identities.
\\
\\
{\sl Keywords:}\quad Pfaffian; Perfect matching; Skew adjacency
matrix; Pl\"ucker relation.
\end{abstract}
\section{Introduction}
\label{sec:intro} \label{sec:Pfaffians} \hspace*{\parindent} Let
$A=(a_{ij})_{n\times n}$ be a skew symmetric matrix of order $n$
and $n$ is even. Suppose that
$\pi=\{(s_1,t_1),(s_2,t_2),\ldots,(s_{\frac{n}{2}},t_{\frac{n}{2}})\}$
is a partition of $[n]$, that is, $[n]=\{s_1,t_1\}\cup
\{s_2,t_2\}\cup \ldots\cup \{s_{\frac{n}{2}},t_{\frac{n}{2}}\}$,
where $[n]=\{1,2,\ldots,n\}$. Define:
$$b_{\pi}=sgn(s_1t_1s_2t_2\ldots s_{\frac{n}{2}}t_{\frac{n}{2}})
\prod_{l=1}^{\frac{n}{2}}a_{s_lt_l},$$ where
$sgn(s_1t_1s_2t_2\ldots s_{\frac{n}{2}}t_{\frac{n}{2}})$ denotes
the sign of the permutation $s_1t_1s_2t_2\ldots
s_{\frac{n}{2}}t_{\frac{n}{2}}$. Note that $b_{\pi}$ depends
neither on the order in which the classes of the partition are
listed nor on the order of the two elements of a class. So
$b_{\pi}$ indeed depends only on the choice of the partition
$\pi$. The Pfaffian of $A$, denoted by $Pf(A)$, is defined as
$$
Pf(A)=\sum_{\pi}b_{\pi},$$ where the summation is over all
partitions of $[n]$, which are of the form of $\pi$. For the sake
of convenience, we define the Pfaffian of $A$ to be zero if $A$ is
a skew symmetric matrix of odd order. The following result is well
known:
\begin{prop}[Cayley Theorem, \cite{Cayley1848}]
For any skew symmetric matrix $A=(a_{ij})_{n\times n}$ of order
$n$, we have
$$\det(A)=[Pf(A)]^2.$$
\end{prop}

Suppose that $G=(V(G),E(G))$ is a weighted graph with the vertex
set $V(G)=\{1,2,\ldots,n\}$, the edge set $E(G)=\{e_1,e_2,\ldots,
e_m\}$ and the edge$-$weight function $\omega: E(G)\longrightarrow
\mathcal R$, where $\omega(e):=\omega_e=a_{ij}$ ($\neq 0$) if
$e=(i,j)$ is an edge of $G$ and $\omega_e=a_{ij}=0$ otherwise, and
$\mathcal R$ is the
set of real numbers. 
Suppose $G^e$ is an orientation of $G$. Let
$A(G^e)=(b_{ij})_{n\times n}$ be the matrix of order $n$ defined
as follows:
$$b_{ij}=\left\{
\begin{array}{cc}
a_{ij} & \mbox{if}\  (i,j) \ \mbox{is an arc in} \ G^e,\\
-a_{ij} & \mbox{if}\ (j,i)\ \mbox{is an arc in}\ G^e,\\
0 & \ \mbox{otherwise}.
\end{array}
\right.
$$
$A(G^e)$ is called the skew adjacency matrix of $G^e$ (see
\cite{Lova86}). Obviously, $A(G^e)$ is a skew symmetric matrix,
that is, $(A(G^e))^T=-A(G^e)$.

Given a skew symmetric matrix $A=(a_{ij})_{n\times n}$ with $n$
even, let $G=(V(G),E(G))$ be a weighted graph with the vertex set
$V(G)=\{1,2,\ldots,n\}$, where $e=(i,j)$ is an edge of $G$ if and
only if $a_{ij}\neq 0$, and the edge$-$weight function is defined
as $\omega_e=|a_{ij}|$ if $e=(i,j)$ is an edge of $G$ and
$\omega_e=0$ otherwise. Define $G^e$ as the orientation of $G$ in
which the direction of every edge $e=(i,j)$ of $G$ is from
vertices $i$ to $j$ if $a_{ij}>0$ and from vertices $j$ to $i$
otherwise. We call $G^e$ to be the corresponding directed graph of
$A$. Obviously, $A=A(G^e)$. It is not difficult to see that the
Pfaffian $Pf(A)$ of $A$ can be defined as
$$Pf(A)=\sum_{\pi\in \mathcal M(G)}b_{\pi},$$
where the summation is over all perfect matchings
$\pi=\{(s_1,t_1),(s_2,t_2),\ldots,(s_{\frac{n}{2}},t_{\frac{n}{2}})\}$
of $G$, and $b_{\pi}$ is the product of all $\omega_{(s_i,t_i)}$
for $1\leq i\leq \frac{n}{2}$.

Pfaffians have been studied for almost two hundred years (see
\cite{Knuth96,Stem90} for a history), and continue to find
numerous applications, for example in matching theory
\cite{Lova86} and in the enumeration of plane partitions
\cite{Stem90}. It is interesting to extend Leclerc's combinatorics
of relations for determinants \cite{Lecl93} to the analogous rules
for Pfaffians. By tools from multilinear algebra Dress and Wenzel
\cite{DW95} gave an elegant proof of an identity concerning
pfaffians of skew symmetric matrices, which yields the
Grassmann-Pl\"ucker identities (for more details see
\cite{Wenz91}, Sect. 7). Okada \cite{Okad05} presented a Pfaffian
identity involving elliptic functions, whose rational limit gives
a generalization of Schur's Pfaffian identity. Knuth
\cite{Knuth96} used a combinatorial method to give an elegant
proof of a classical Pfaffian identity found in \cite{Tann78}.
Hamel \cite{Hame01} followed Knuth's approach and introduced other
combinatorial methods to prove a host of Pfaffian identities from
physics in \cite{Hiro89,Ohta93,TH96}. Hamel also provided a
combinatorial proof of a result in \cite{Srin94} and a new
vector-based Pfaffian identity and gave an application to the
theory of symmetric functions by proving an identity for Schur
$Q$-functions. For some related recent results see also
\cite{IW95,IW05,LT02,Okad98}.

This paper is inspired by two results, one of which is that we can
use the Pfaffian method to enumerate perfect matchings of plane
graphs (see \cite{Kast63,Kast67}). Inspired by the Dodgson's
Determinant-Evaluation Rule in \cite{Dodg66} and the Pl\"ucker
relations for Pfaffians, Propp \cite{Prop03}, Kuo \cite{Kuo04} and
Yan et al \cite{YZ05} obtained a method of graphical
vertex-condensation for enumerating perfect matchings of plane
bipartite graphs. The second is that by using the Matching
Factorization Theorem in \cite{Ciuc97} Yan et al \cite{YYZ05}
found a method of graphical edge-condensation for counting perfect
matchings of plane graphs. It is natural to ask whether there
exist some Pfaffian identities completely different from the
Pl\"ucker relations, which can result in some formulas for the
method of graphical edge-condensation for enumerating perfect
matchings of plane graphs. The results in Section 3 answer this
question in the affirmative. As applications, we obtain two new
determinant identities in Section 4.1 and we prove a quadratic
relation for the number of perfect matchings of plane graphs in
Section 4.2, which has a simpler form than the formula in
\cite{YYZ05}.
\section{Some Lemmas}
\hspace*{\parindent} In order to present the following lemmas, we
need to introduce some notation and terminology. If $I$ is a
subset of ${\bf [n]}$, we use $A_{I}$ to denote the minor of $A$
by deleting rows and columns indexed by $I$. If
$I=\{i_1,i_2,\ldots,i_l\}\subseteq {\bf [n]}$ and
$i_1<i_2<\ldots<i_l$, we use $Pf_A(i_1i_2\ldots i_l)=: Pf_{A}(I)$
to denote the Pfaffian of $A_{{\bf [n]}\backslash I}$. Following
Knuth's notation in \cite{Knuth96}, for two words $\alpha$ and
$\beta$ we define $s(\alpha,\beta)$ to be zero if either $\alpha$
or $\beta$ has a repeated letter, or if $\beta$ contains a letter
not in $\alpha$. Or, if these are not the case, $s(\alpha,\beta)$
denotes the sign of the permutation that takes $\alpha$ into the
word $\beta(\alpha\backslash \beta)$ (where $\alpha\backslash
\beta$ denotes the word that remains when the elements of $\beta$
are removed from $\alpha$). Let $S$ be a subset of
$\{1,2,\ldots,n\}$. We call $S$ an even subset if $|S|$ is even
and an odd one otherwise.

Dress et al \cite{DW95} used tools from multilinear algebra to
prove a Pfaffian identity, which was found by Wenzel
\cite{Wenz91}, as follows:
\begin{lem}[Wenzel \cite{Wenz91} and Dress et al \cite{DW95}]
For any two subsets $I_1, I_2\subseteq [n]$ of odd cardinality and
elements $i_1,i_2,\ldots,i_t\in [n]$ with $i_1<i_2\ldots<i_t$ and
$\{i_1,i_2,\ldots,i_t\}=I_1\triangle I_2=: (I_1\backslash I_2)\cup
(I_2\backslash I_1)$, if $A=(a_{ij})_{n\times n}$ is a skew
symmetric matrix with $n$ even, then
\begin{center}
$\sum\limits_{\tau=1}^t(-1)^{\tau}Pf_A(I_1\triangle
\{i_{\tau}\})Pf(I_2\triangle \{i_{\tau}\})=0.$
\end{center}
\end{lem}

A direct result of Lemma 2.1 is the following lemma, which will
play an important role in the proofs of our main results.
\begin{lem}
Suppose that $A=(a_{ij})_{n\times n}$ is a skew symmetric matrix
with $n$ even and $\alpha$ is an even subset of $\bf [n]$. Let
$\beta=\{i_1,i_2,\ldots,i_{2p}\} \subseteq [n]\backslash \alpha$,
where $i_1<i_2<\ldots<i_{2p}$. Then, for any fixed $s\in [2p]$, we
have
\begin{center}
$Pf_A(\alpha)Pf_A(\alpha\beta)=\sum\limits_{l=1}^{2p}(-1)^{l+s+1}Pf_A(\alpha
i_si_l)Pf_A(\alpha\beta\backslash i_si_l),$
\end{center}
where $Pf_A(\alpha i_si_s)=0$.
\end{lem}

The following result is a special case of Lemma 2.2.
\begin{cor}
Suppose that $A=(a_{ij})_{n\times n}$ is a skew symmetric matrix
and $\{i,j,k,l\}\subseteq {\bf [n]}$. Then
$$Pf(A_{\{i,j,k,l\}})Pf(A)=Pf(A_{\{i,j\}})Pf(A_{\{k,l\}})-
Pf(A_{\{i,k\}})Pf(A_{\{j,l\}})+Pf(A_{\{i,l\}})Pf(A_{\{j,k\}}).\eqno{(2.1)}$$
\end{cor}
\begin{rem}
There exists a similar formula on the determinant to Corollary 2.1
as follows, which is called the Dodgson's Determinant$-$Evaluation
Rule (see \cite{Dodg66}):
$$\det(A_{\{1,n\}})\det(A)=\det(A_{11})\det(A_{nn})-\det(A_{1n})\det(A_{n1}), \eqno{(2.2)}$$
where $A$ is an arbitrary matrix of order $n$ and $A_{ij}$ is the
minor of $A$ by deleting the $i$-th row and the $j$-th column.
\end{rem}

The following result shows the relation between the Pfaffian and
the determinant.
\begin{lem}[Godsil \cite{Gods93}]
Let $A$ be a square matrix of order $n$. Then
$$Pf\left(\begin{array}{cc}
0 & A\\
-A^T & 0
\end{array}
\right)=(-1)^{\frac{1}{2}n(n-1)}\det (A).$$
\end{lem}

Let $A=(a_{st})_{n\times n}$ be a skew symmetric matrix of order
$n$ and $G^e$ the corresponding directed graph. Suppose $(i,j)$ is
an arc in $G^e$ and hence $a_{ij}>0$. Let $\overline G^e$ be a
directed graph with vertex set $\{1,2,\ldots, n+1,n+2\}$ obtained
from $G^e$ by deleting the arc $(i,j)$ and adding three arcs
$(i,n+1), (n+1,n+2)$ and $(n+2,j)$ with weights $\sqrt {a_{ij}}, 1
$ and $\sqrt{a_{ij}}$, respectively (see Figures 1(a) and (b) for
an illustration). For convenience, if $a_{ij}=0$ we also regard
$\overline G^e$ as a directed graph obtained from $G^e$ by adding
three arcs $(i,n+1), (n+1,n+2)$ and $(n+2,j)$ with weights $0, 1 $
and $0$. The following lemma will play a key role in the proofs of
our main results.
\begin{figure}[htbp]
  \centering
  \scalebox{0.75}{\includegraphics{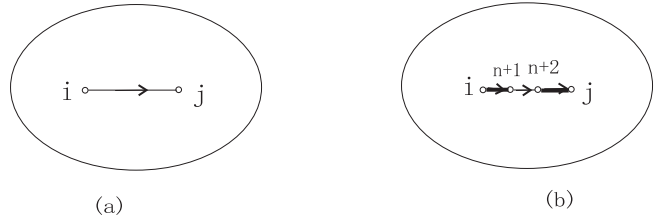}}
  \caption{\ (a)\ The directed graph $G^e$.
  \ (b)\ The directed graph $\overline G^e$.}
\end{figure}
\begin{lem}
Suppose that $A=(a_{st})_{n\times n}$ is a skew symmetric matrix
and $G^e$ is the corresponding directed graph. Let $\overline G^e$
be the directed graph with $n+2$ vertices defined above and
$A(\overline G^e)$ the skew adjacency matrix of $\overline G^e$.
Then
$$Pf(A)=Pf(A(\overline G^e)).$$
\end{lem}

{\bf Proof}\ \ Let $G$ and $\overline G$ be the underlying graphs
of $G^e$ and $\overline G^e$, and let $A(G^e)$ be the skew
adjacency matrices of $G^e$. Hence $A(G^e)=(a_{st})_{n\times n}$
and $A(\overline G^e)=(b_{st})_{(n+2)\times (n+2)}$, where
$$
b_{st}=\left\{\begin{array}{ll}
a_{st} & \ \ \mbox{if}\ 1\leq s,t\leq n \ \mbox{and}\ (s,t)\neq (i,j), (j,i),\\
\sqrt {a_{ij}} & \ \ \mbox{if}\ (s,t)=(i,n+1)\  \mbox{or} \ (n+2,j),\\
-\sqrt {a_{ij}} & \ \ \mbox{if}\ (s,t)=(n+1,i)\  \mbox{or} \ (j,n+2),\\
1 & \ \ \mbox{if}\ (s,t)=(n+1,n+2),\\
-1 & \ \ \mbox{if}\ (s,t)=(n+2,n+1),\\
0 & \ \ \mbox{otherwise}.
\end{array}
\right.
$$

By the definitions above, we have
$$Pf(A)=Pf(A(G^e)).$$
Hence we only need to prove $$Pf(A(G^e))=Pf(A(\overline G^e)).$$
Note that, by the definition of the Pfaffian, we have

$$Pf(A(G^e))=\sum\limits_{\pi\in \mathcal M(G)}b_{\pi},\ \
Pf(A(\overline G^e))= \sum\limits_{\overline \pi\in \mathcal
M(\overline G)}b_{\overline \pi},$$ where $\mathcal M(G)$ and
$\mathcal M(\overline G)$ denote the sets of perfect matchings of
$G$ and $\overline G$.

We partition the sets of perfect matchings of $G$ and $\overline
G^e$ as follows:
$$\mathcal M(G)=\mathcal M_1\cup \mathcal M_2,\ \
\mathcal M(\overline G)=\overline {\mathcal M_1}\cup \overline
{\mathcal M_2},$$ where $\mathcal M_1$ is the set of perfect
matchings of $G$ each of which contains edge $e=(i,j)$, $\mathcal
M_2$ is the set of perfect matchings of $G$ each of which does not
contain edge $e=(i,j)$, $\overline {\mathcal M_1}$ is the set of
perfect matchings of $\overline G$ each of which contains both of
edges $(i,n+1)$ and $(n+2,j)$, and $\overline {\mathcal M_2}$ is
the set of perfect matchings of $\overline G$ each of which
contains edge $(n+1,n+2)$.

Suppose $\pi$ is a perfect matching of $G$. If $\pi\in \mathcal
M_1$, then there exists uniquely a perfect matching $\pi'$ of
$G-i-j$ such that $\pi=\pi'\cup \{(i,j)\}$. It is clear that there
is a natural way to regard $\pi'$ as a matching of $\overline G$.
Define: $\overline \pi=\pi'\cup \{(i,n+1),(n+2,j)\}$. Hence
$\overline \pi\in \overline {\mathcal M_1}$. Similarly, if $\pi\in
\mathcal M_2$, we can define: $\overline \pi=\pi\cup
\{(n+1,n+2)\}$ and hence $\overline {\pi}\in \overline {\mathcal
M_2}$. It is not difficult to see that the mapping $f: \pi
\longmapsto \overline \pi$ between $\mathcal M(G)$ and $\mathcal
M(\overline G)$ is bijective.

Hence we only need to prove that for any perfect matching $\pi$ of
$G$ we have $b_{\pi}=b_{\overline {\pi}}$. By the definition of
$\overline \pi$, if
$\pi=\{(s_1,t_1),(s_2,t_2),\ldots,(s_{l-1},t_{l-1}), (i,j),
(s_{l+1},t_{l+1}),\ldots,(s_{\frac{n}{2}},t_{\frac{n}{2}})\}\in
\mathcal M_1$, then $\overline \pi=
\{(s_1,t_1),(s_2,t_2),\ldots,(s_{l-1},t_{l-1}),
(i,n+1),(n+2,j),(s_{l+1},t_{l+1}),\ldots,(s_{\frac{n}{2}},t_{\frac{n}{2}})\}\in
\overline{\mathcal M_1}$. Note that
$$
sgn(s_1t_1\ldots s_{l-1}t_{l-1}ijs_{l+1}t_{l+1}\ldots
s_{\frac{n}{2}}t_{\frac{n}{2}})= sgn(s_1t_1\ldots
s_{l-1}t_{l-1}i(n+1)(n+2)js_{l+1}t_{l+1}\ldots
s_{\frac{n}{2}}t_{\frac{n}{2}}),
$$
$b_{s_1t_1}\ldots
b_{s_{l-1}t_{l-1}}b_{i(n+1)}b_{(n+2)j}b_{s_{l+1}t_{l+1}}\ldots
b_{s_{\frac{n}{2}}b_{\frac{n}{2}}}= $
$$ a_{s_1t_1}\ldots a_{s_{l-1}t_{l-1}}\sqrt
{a_{ij}}\sqrt{a_{ij}}a_{s_{l+1}t_{l+1}}\cdots
a_{s_{\frac{n}{2}}t_{\frac{n}{2}}}=a_{s_1t_1}\ldots
a_{s_{l-1}t_{l-1}}a_{ij}a_{s_{l+1}t_{l+1}}\ldots
a_{s_{\frac{n}{2}}t_{\frac{n}{2}}}.
$$ Thus we have showed that if
$\pi\in \mathcal M_1$ then we have $b_{\pi}=b_{\overline \pi}$.
Similarly, we can prove that if $\pi\in \mathcal M_2$ then we have
$b_{\pi}=b_{\overline \pi}$. So we have proved that
$Pf(A(G^e))=Pf(A(\overline G^e))$, and the lemma follows. $\hfill
\blacksquare$

\section{New Pfaffian identities}
\hspace*{\parindent}
We first need to introduce some notation. In this section, we
assume that $A=(a_{ij})_{n\times n}$ is a skew symmetric matrix
with $n$ even. Suppose $E=\{(i_l,j_l)|l=1,2,\ldots,k\}$ is a
subset of ${\bf [n]\times [n]}$ such that $i_1\leq i_2\leq \ldots
\leq i_k$ and $i_l<j_l$ for $1\leq l\leq k$. We define a new skew
symmetric matrix $E(A)$ of order $n$ from $A$ and $E$ as follows:
$$
E(A)=(b_{ij})_{n\times n},\ \ b_{ij}=\left\{
\begin{array}{ll}
a_{ij} & \ \mbox{if}\ (i,j)\notin E\  \mbox {and}\ i<j,\\
-a_{ji} & \ \mbox{if}\ (j,i)\notin E \ \mbox{and}\ i>j,\\
0 &\ \mbox{otherwise}.
\end{array}
\right.
$$
By the definition of $E(A)$, it is obtained from $A$ by replacing
all $(i_l,j_l)$ and $(j_l,i_l)-$entries with $zeros$ and not
changing the other entries and hence it is a skew symmetric
matrix. For example, if $A=(a_{ij})_{4\times 4}$ is a skew
symmetric matrix and $E=\{(1,4),(2,3),(3,4)\}$, then
$$E(A)=\left(\begin{array}{cccc}
0 & a_{12} & a_{13} & 0\\
-a_{12} & 0 & 0 & a_{24}\\
-a_{13} & 0 & 0 & 0\\
0 & -a_{24} & 0 & 0
\end{array}
\right).
$$

Now, we can state one of our main results as follows.
\begin{thm}
Suppose $A=(a_{ij})_{n\times n}$ is a skew symmetric matrix of
order $n$ and $E=\{(i_l,j_l)|l=1,2,\ldots,k\}$ is a non empty
subset of ${\bf [n]\times [n]}$ such that $i_1\leq i_2\leq \ldots
\leq i_k,\ i_l<j_l$ for $l\in [k]$. Then, for any fixed $p \in
[k]$, we have
$$
Pf(E(A))Pf(A)=Pf(E_p(A))Pf(\overline{E_p}(A))+
$$
$$a_{i_pj_p}\sum_{1\leq l\leq k, l\neq p}a_{i_lj_l}\left[f(p,l)
Pf(E(A)_{\{i_p,j_l\}})Pf(A_{\{j_p,i_l\}})
-g(p,l)Pf(E(A)_{\{i_p,i_l\}})Pf(A_{\{j_p,j_l\}})\right],$$
where $E_p=E\backslash \{(i_p,j_p)\},
\overline{E_p}=\{(i_p,j_p)\}, f(p,l)=s([n],i_pj_l)s([n],j_pi_l)$
and $g(p,l)=s([n],i_pi_l)s([n],j_pj_l)$.
\end{thm}

{\bf Proof}\ \ Let $G^e$ be the corresponding directed graph of
$A$ defined as above, whose vertex set is ${\bf [n]}$. Let
$\overline G^e$ be the directed graph with the vertex set $[n+2k]$
obtained from $G^e$ by replacing each arc between every pair of
vertices $i_l$ and $j_l$ with three arcs $(i_l,n+2l-1),
(n+2l-1,n+2l)$ and $(n+2l,j_l)$ with weights $\sqrt {a_{i_lj_l}},
1$ and $\sqrt{a_{i_lj_l}}$ if $(i_l,j_l)$ is an arc of $G^e$ and
with three arcs $(j_l, n+2l-1), (n+2l-1, n+2l)$ and $(n+2l,i_l)$
with weights $\sqrt {a_{j_li_l}}, 1$ and $\sqrt{a_{j_li_l}}$ if
$(j_l,i_l)$ is an arc of $G^e$, respectively.  For the case
$a_{i_lj_l}>0$ for $1\leq l\leq k$, Figure 2 (a) and (b)
illustrate the procedure constructing $\overline G^e$ from $G^e$.
Suppose $\overline A=A(\overline G^e)$ is the skew adjacency
matrix of $\overline G^e$.
\begin{figure}[htbp]
  \centering
  \scalebox{0.8}{\includegraphics{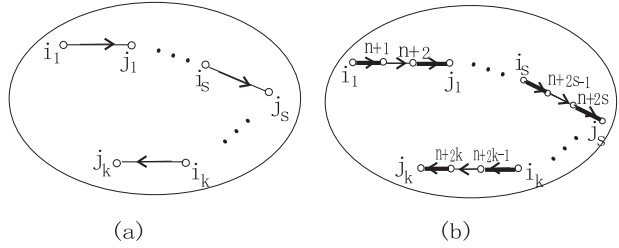}}
  \caption{\ (a)\ The directed graph $G^e$.
  \ (b)\ The directed graph $\overline G^e$.}
\end{figure}

Take $\alpha=[n],
\beta=\{n+1,n+2,\ldots,n+2k\}=\{x_i|x_i=n+i,1\leq i\leq 2k\}$.
Take $q=2p-1$. Hence $x_q=n+2p-1$ and $(-1)^{l+q+1}=(-1)^l$. By
Lemma 2.2, we have
$$Pf_{\overline A}(\alpha)Pf_{\overline A}(\alpha\beta)=
\sum_{l=1}^{2k}(-1)^lPf_{\overline A}(\alpha x_qx_l) Pf_{\overline
A}(\alpha\beta\backslash x_qx_l). \eqno{(3.1)}
$$
By the definitions of $E(A)$ and $G^e$ and Lemma 2.4, we have
$$Pf_{\overline A}(\alpha)=Pf(E(A)), \ \ Pf_{\overline A}(\alpha\beta)=Pf(A).
\eqno{(3.2)}
$$

We set
\begin{center}
$a_{l'}=-Pf_{\overline A}(\alpha x_qx_{2l'-1})Pf_{\overline
A}(\alpha\beta\backslash x_qx_{2l'-1}),$\\
$b_{l'}=Pf_{\overline A}(\alpha x_qx_{2l'})Pf_{\overline
A}(\alpha\beta\backslash x_qx_{2l'}),$
\end{center}
that is,
$$a_{l'}=-Pf_{\overline A}(\alpha (n+2p-1)(n+2l'-1))
Pf_{\overline A}(\alpha\beta\backslash
(n+2p-1)(n+2l'-1)),\eqno{(3.3)}$$
$$b_{l'}=Pf_{\overline A}(\alpha (n+2p-1)(n+2l'))
Pf_{\overline A}(\alpha\beta\backslash (n+2p-1)(n+2l')).
\eqno{(3.4)}
$$
By Lemma 2.4, it is not difficult to see that
$$b_p=Pf_{\overline A}(\alpha x_qx_p)Pf_{\overline A}(\alpha\beta\backslash
\{x_qx_p\})=Pf(E_p(A))Pf(\overline{E_p}(A)). \eqno{(3.5)}$$ Note
that $a_p=0$. Hence we have
$$Pf_{\overline A}(\alpha)Pf_{\overline A}(\alpha\beta)=\sum\limits_{l=1}^{2k}(-1)^{l}Pf_{\overline A}
(\alpha x_qx_l)Pf_{\overline A}(\alpha\beta\backslash x_qx_l)$$$$
=Pf(E_p(A))Pf(\overline{E_p}(A))+\sum_{1\leq l'\leq k,l'\neq
p}(a_{l'}+b_{l'}).\eqno{(3.6)}$$

Obviously, if $a_{i_pj_p}=0$ then theorem is trivial. Hence we may
assume that $a_{i_pj_p}\neq 0$.

First, we prove that if $a_{i_pj_p}>0$ then the theorem holds.
From $(3.2)$ and $(3.6)$ it suffices to prove the following claim:

{\bf Claim}\ \ For any $l'\in [k]$ and $l'\neq p$, if
$a_{i_pj_p}>0$ we have
$$a_{l'}+b_{l'}=
s(([n]),i_pj_{l'})s([n],j_pi_{l'})a_{i_pj_p}a_{i_{l'}j_{l'}}
Pf(E(A)_{\{i_p,j_{l'}\}})Pf(A_{\{j_p,i_{l'}\}})-$$
$$s([n],i_pi_{l'})s([n],j_pj_{l'})
a_{i_pj_p}a_{i_{l'}j_{l'}}Pf(E(A)_{\{i_p,i_{l'}\}})
Pf(A_{\{j_p,j_{l'}\}}).\eqno{(3.7)}
$$

Suppose $a_{i_pj_p}>0$. Then $(i_p,n+2p-1), (n+2p-1,n+2p)$ and
$(n+2p,j_p)$ are three arcs of $\overline G^e$ with weights
$\sqrt{a_{i_pj_p}}, 1$ and $\sqrt{a_{i_pj_p}}$. We need to
consider two cases:

$(a)$\ \ $a_{i_{l'}j_{l'}}\geq 0$;

$(b)$\ \ $a_{i_{l'}j_{l'}}<0$.

If $a_{i_{l'}j_{l'}}\geq 0$, then $(i_{l'},n+2l'-1),
(n+2l'-1,n+2l')$ and $(n+2l',j_{l'})$ are three arcs of $\overline
G^e$ with weights $\sqrt{a_{i_{l'}j_{l'}}}, 1$ and
$\sqrt{a_{i_{l'}j_{l'}}}$. Suppose $X$ is a subset of the vertex
set of $\overline{G}$. Let $\overline{G}[X]^e$ be the directed
subgraph of $\overline{G}^e$ induced by $X$ and $\overline{G}[X]$
the underlying graph of $\overline{G}[X]^e$. Note that
$\overline{G}[\alpha(n+2l'-1)(n+2p-1)]$ contains two pendant edges
$(i_{l'},n+2i_{l'}-1)$ and $(i_p,n+2p-1)$. Each perfect matching
$\pi$ of $\overline{G}[\alpha(n+2l'-1)(n+2p-1)]$ can be denoted by
$\pi=\pi'\cup\{(i_{l'},n+2i_{l'}-1),(i_p,n+2p-1)\}$, where $\pi'$
is a perfect matching of $\overline{G}[\alpha\backslash
\{i_{l'},i_p\}]$. Set

$Pf(A(\overline G[\alpha(n+2l'-1)(n+2p-1)]^e))=\sum\limits_{\pi\in
\mathcal M(\overline G[\alpha(n+2l'-1)(n+2p-1)])}b_{\pi},$

$Pf(A(\overline G[\alpha \backslash
\{i_{l'},i_p\}]^e))=\sum\limits_{\pi'\in \mathcal M(\overline
G[\alpha\backslash \{i_{l'},i_p\}])}b_{\pi'},$\\
where $\mathcal M(G)$ is the set of perfect matchings of a graph
$G$. By the definitions of $b_{\pi}$ and $b_{\pi'}$, it is not
difficult to see that
$$b_{\pi}=sgn(p-l')s([n],i_pi_{l'})\sqrt{a_{i_pj_p}a_{i_{l'}j_{l'}}}b_{\pi'},$$
where $sgn(x)$ denotes the sign of $x$. By the definition of
$E(A)$, we have $$Pf(A(\overline G[\alpha \backslash
\{i_{l'},i_p\}]^e))=Pf(E(A)_{\{i_{l'},i_p\}}).$$ Hence we have
proved the following:
$$Pf_{\overline A}(\alpha (n+2l'-1)(n+2p-1))=sgn(p-l')s([n],i_pi_{l'})
\sqrt{a_{i_pj_p}a_{i_{l'}j_{l'}}} Pf(E(A)_{\{i_p,i_{l'}\}}).
\eqno{(3.8)}$$
Similarly, we can prove the following:
$$Pf_{\overline A}(\alpha\beta\backslash (n+2l'-1)(n+2p-1))
=sgn(p-l')s([n],j_pj_{l'})\sqrt{a_{i_pj_p}a_{i_{l'}j_{l'}}}
Pf(A_{\{j_p,j_{l'}\}}); \eqno{(3.9)}$$
$$Pf_{\overline A}(\alpha (n+2l')(n+2p-1))=sgn(l'-p)
s([n],i_pj_{l'})\sqrt{a_{i_pj_p}a_{i_{l'}j_{l'}}}
Pf(E(A)_{\{i_p,j_{l'}\}}); \eqno{(3.10)}$$
$$Pf_{\overline A}(\alpha\beta\backslash (n+2l')(n+2p-1))
=sgn(l'-p)s([n],j_pi_{l'})\sqrt{a_{i_pj_p}a_{i_{l'}j_{l'}}}Pf(A_{\{j_p,i_{l'}\}}).
\eqno{(3.11)}$$ Then $(3.7)$ is immediate from $(3.3),(3.4),
(3.8)-(3.11)$. Hence if $a_{i_{l'}j_{l'}}\geq 0$ then the claim
follows.

If $a_{i_{l'}j_{l'}}<0$, then $(j_{l'},n+2l'-1), (n+2l'-1,n+2l')$
and $(n+2l',i_{l'})$ are three arcs of $\overline G^e$ with
weights $\sqrt{-a_{i_{l'}j_{l'}}}, 1$ and
$\sqrt{-a_{i_{l'}j_{l'}}}$. Similarly, we can prove the following:
$$Pf_{\overline A}(\alpha (n+2l'-1)(n+2p-1))=sgn(p-l')
s([n],i_pj_{l'})\sqrt{a_{i_pj_p}a_{j_{l'}i_{l'}}}
Pf(A(E(A)_{\{i_p,j_{l'}\}}); \eqno{(3.12)}$$
$$Pf_{\overline A}(\alpha\beta\backslash (n+2l'-1)(n+2p-1))
=sgn(p-l')s([n],j_pi_{l'})\sqrt{a_{i_pj_p}a_{j_{l'}i_{l'}}}
Pf(A_{\{j_p,i_{l'}\}}); \eqno{(3.13)}$$
$$Pf_{\overline A}(\alpha (n+2l')(n+2p-1))=sgn(l'-p)
s([n],i_pi_{l'})\sqrt{a_{i_pj_p}a_{j_{l'}i_{l'}}}
Pf(E(A)_{\{i_p,i_{l'}\}}); \eqno{(3.14)}$$
$$Pf_{\overline A}(\alpha\beta\backslash (n+2l')(n+2p-1))
=sgn(l'-p)s([n],j_pj_{l'})\sqrt{a_{i_pj_p}a_{j_{l'}i_{l'}}}
Pf(A_{\{j_p,j_{l'}\}}). \eqno{(3.15)}$$ Then $(3.7)$ is immediate
from $(3.3),(3.4), (3.12)-(3.15)$. Hence if $a_{i_{l'}j_{l'}}<0$
then the claim follows.

Hence we have proved that if $a_{i_pj_p}>0$ then the theorem
holds.

If $a_{i_pj_p}<0$, we consider $Pf(-A)$ and $Pf(-E(A))$. Note that
$(-A)_{i_pj_p}>0$. The result proved above implies that
$$
Pf(-E(A))Pf(-A)=Pf(E_p(-A))Pf(\overline{E_p}(-A))-
a_{i_pj_p}\sum_{1\leq l\leq k, l\neq
p}(-a_{i_lj_l})\times$$$$\left[f(p,l)\times
Pf(E(-A)_{\{i_p,j_l\}})Pf((-A)_{\{j_p,i_l\}})
-g(p,l)Pf(E(-A)_{\{i_p,i_l\}})Pf((-A)_{\{j_p,j_l\}})\right].\eqno{(3.16)}$$
Note that by the definition of the Pfaffian we have
$Pf(-A)=(-1)^{\frac{n}{2}}Pf(A)$. By $(3.16)$, we can show that we
have
$$
Pf(E(A))Pf(A)=Pf(E_p(A))Pf(\overline{E_p}(A))+
$$
$$
a_{i_pj_p}\sum_{1\leq l\leq k, l\neq
p}a_{i_lj_l}\left[f(p,l)Pf(E(A)_{\{i_p,j_l\}})Pf(A_{\{j_p,i_l\}})
-g(p,l)Pf(E(A)_{\{i_p,i_l\}})Pf(A_{\{j_p,j_l\}})\right],
$$
which implies that if $a_{i_pj_p}<0$ then the theorem also holds.

Hence we have proved the theorem. $\hfill \blacksquare$
\begin{cor}
With the same notation as Theorem 3.1, for any fixed $p \in [k]$,
$$
Pf(E(A))Pf(A)=Pf(E_p(A))Pf(\overline{E_p}(A))+
$$
$$a_{i_pj_p}\sum_{1\leq l\leq k, l\neq p}a_{i_lj_l}\left[f(p,l)
Pf(E(A)_{\{j_p,i_l\}})Pf(A_{\{i_p,j_l\}})
-g(p,l)Pf(E(A)_{\{j_p,j_l\}})Pf(A_{\{i_p,i_l\}})\right].$$
\end{cor}

{\bf Proof}\ \ Let $A^T$ be the transpose of $A$. Note that
$Pf(A^T)=(-1)^{\frac{n}{2}}Pf(A)$. The corollary follows
immediately from Theorem 3.1 by considering the transpose of $A$.
$\hfill \blacksquare$\\

The following result is a special case of Theorem 3.1 and
Corollary 3.2.
\begin{cor}
Suppose $A=(a_{ij})_{n\times n}$ is a skew symmetric matrix of
order $n$ and $E=\{(i_l,j_l)|l=1,2,\ldots,k\}$ is a non empty
subset of ${\bf [n]\times [n]}$ such that $i_1<j_1<i_2<j_2<\ldots
i_l<j_l<\ldots<i_k<j_k$. Then

$Pf(E(A))Pf(A)-Pf(E_1(A))Pf(\overline{E_1}(A))$

$=a_{i_1j_1}\sum\limits_{l=2}^ka_{i_lj_l}\left[
Pf(E(A)_{\{i_1,j_l\}})Pf(A_{\{j_1,i_l\}})
-Pf(E(A)_{\{i_1,i_l\}})Pf(A_{\{j_1,j_l\}})\right]$

$=a_{i_1j_1}\sum\limits_{l=2}^ka_{i_lj_l}\left[
Pf(E(A)_{\{j_1,i_l\}})Pf(A_{\{i_1,j_l\}})
-Pf(E(A)_{\{j_1,j_l\}})Pf(A_{\{i_1,i_l\}})\right].$
\end{cor}
\begin{rem}
The Pfaffian identities in Theorem 3.1 and Corollaries 3.2 and 3.3
express the product of Pfaffians of two skew symmetric matrices
$E(A)$ and $A$ in terms of the Pfaffians of the minors of $E(A)$
and $A$, where $E(A)$ is a skew symmetric matrix obtained from $A$
by replacing some non zero entries $a_{i_lj_l}$ and $a_{j_li_l}$
of $A$ with zeros. On the other hand, an obvious observation in
the Pfaffian identities known before, which belong to the
Pl\"ucker relations, is that the related matrices are either a
skew symmetric matric $A$ or some minors of $A$. Hence the
Pfaffian identities in Theorem 3.1 and Corollaries 3.2 and 3.3 are
completely new and different from the Pl\"ucker relations.
\end{rem}
\begin{ex}
Let $A=(a_{ij})_{4\times 4}$ and $E=\{(1,2),(3,4)\}$. Then, by
Corollary 3.3, we have

{\small $Pf\left(\begin{array}{cccc}
0 & a_{12} & a_{13} & a_{14}\\
-a_{12} & 0 & a_{23} & a_{24}\\
-a_{13} & -a_{23} & 0 & a_{34}\\
-a_{14} & -a_{24} & -a_{34} & 0
\end{array}
\right)Pf\left(\begin{array}{cccc}
0 & 0 & a_{13} & a_{14}\\
0 & 0 & a_{23} & a_{24}\\
-a_{13} & -a_{23} & 0 & 0\\
-a_{14} & -a_{24} & 0 & 0
\end{array}
\right)=$}

{\small $Pf\left(\begin{array}{cccc}
0 & 0 & a_{13} & a_{14}\\
0 & 0 & a_{23} & a_{24}\\
-a_{13} & -a_{23} & 0 & a_{34}\\
-a_{14} & -a_{24} & -a_{34} & 0
\end{array}
\right)Pf\left(\begin{array}{cccc}
0 & a_{12} & a_{13} & a_{14}\\
-a_{12} & 0 & a_{23} & a_{24}\\
-a_{13} & -a_{23} & 0 & 0\\
-a_{14} & -a_{24} & 0 & 0
\end{array}
\right)+$\\
$a_{12}a_{34}Pf\left(\begin{array}{cc} 0 & a_{23}\\
-a_{23} & 0\end{array}\right)Pf\left(\begin{array}{cc} 0 & a_{14}\\
-a_{14} & 0\end{array}\right)+a_{12}a_{34}Pf\left(\begin{array}{cc} 0 & a_{24}\\
-a_{24} & 0\end{array}\right)Pf\left(\begin{array}{cc} 0 & a_{13}\\
-a_{13} & 0\end{array}\right).$}
\end{ex}
\section{Applications}
\hspace*{\parindent}
As applications of some results in Section 3, we obtain some
determinant identities different from the Pl\"ucker relations in
Section 4.1 and we prove a quadratic relation for the number of
perfect matchings of plane graphs in Section 4.2, which has a
simpler form than the formula in \cite{YYZ05}.
\subsection{New determinant identities}
\hspace*{\parindent}
We first need to introduce some notation and terminology.
Throughout this subsection, we will assume $A=(a_{ij})_{n\times
n}$ is an arbitrary  matrix of order $n$ and $E=\{(i_l,j_l)|1\leq
l\leq k\}\subseteq [n]\times [n]$, where $a_{i_lj_l}\neq 0$.
Define a new matrix of order $n$ from $A$ and $E$, denoted by
$E[A]=(b_{st})_{n\times n}$, where $ b_{st}=\left
\{\begin{array}{ll}
a_{st} & \ \ \mbox{if}\ \ (s,t)\notin E,\\
0 & \ \ \mbox{otherwise}.
\end{array}
\right. $. In other words, $E[A]$ is an $n\times n$ matrix
obtained from $A$ by replacing all entries $a_{i_lj_l}$ for $1\leq
l\leq k$ with zeros and not changing the other entries. For
example, if $A=(a_{ij})_{4\times 4}, E=\{(1,2),(2,2),(3,1)\}$, by
the definition of $E[A]$ we have
$$
E[A]=\left (\begin{array}{cccc}
a_{11} & 0 & a_{13} & a_{14}\\
a_{21} & 0 & a_{23} & a_{24}\\
0 & a_{32} & a_{33} & a_{34}\\
a_{41} & a_{42} & a_{43} & a_{44}
\end{array}
\right), \{(3,4)\}[A]=\left (\begin{array}{cccc}
a_{11} & a_{12} & a_{13} & a_{14}\\
a_{21} & a_{22} & a_{23} & a_{24}\\
a_{31} & a_{32} & a_{33} & 0\\
a_{41} & a_{42} & a_{43} & a_{44}
\end{array}
\right).
$$

Now, we start to prove the following:
\begin{lem}
If $A=(a_{ij})_{n\times n}$ is a matrix of order $n$ and
$A^*=\left(\begin{array}{cc}
0 & A\\
-A^T & 0
\end{array}
\right)
$, then, for any $i,j\in [n]$, $i\neq j$, we have\\
$(i)$\ \ $Pf(A^*_{\{i,j\}})=0,$\\
$(ii)$\ \ $Pf(A^*_{\{i,n+j\}})=(-1)^{\frac{1}{2}(n-1)(n-2)}\det(A_{ij}),$\\
$(iii)$\ \ $Pf(A^*_{\{n+i,j\}})=(-1)^{\frac{1}{2}(n-1)(n-2)}\det(A_{ji}),$\\
where $A_{ij}$ denotes the minor of $A$ obtained by deleting the
$i-$th row and $j-$th column from $A$.
\end{lem}

{\bf Proof}\ \ Note that $A^*_{\{i,j\}}=\left(\begin{array}{cc}
0 & B\\
-B^T & 0
\end{array}
\right)$, where $B$ is an $(n-2)\times n$ matrix obtained from $A$
by deleting two rows indexed by $i$ and $j$. Obviously, $\det
\left(\begin{array}{cc}
0 & B\\
-B^T & 0
\end{array}
\right)=0$. Hence by Cayley Theorem we have
$\left[Pf(A^*_{\{i,j\}})\right]^2=\det \left(\begin{array}{cc}
0 & B\\
-B^T & 0
\end{array}
\right)=0$, which implies that $Pf(A^*_{\{i,j\}})=0$. Similarly,
by Lemma 2.3 we can prove $(ii)$ and $(iii)$. Hence the lemma
follows.\ \ \ \ $\hfill \blacksquare$

\begin{thm}
Let $A=(a_{ij})_{n\times n}$ be a matrix of order $n$ and
$E=\{(i_l,j_l)|1\leq l\leq k\}$ a non empty subset of $[n]\times
[n]$, where $i_1\leq i_2\leq\ldots\leq i_k$. Then for a fixed
$p\in [k]$ we have
$$ \det(E[A])\det(A)$$
$$=\det(E_p[A]) \det(\overline{E_p}[A])-
\sum\limits_{1\leq l\leq k,l\neq p}(-1)^{i_l+j_l+i_p+j_p}
a_{i_pj_p}a_{i_lj_l}\det(E[A]_{i_pj_l})\det(A_{i_lj_p}), $$ where
$E_p=E\backslash\{(i_p,j_p)\}$ and $\overline{E_p}=\{(i_p,j_p)\}$.
\end{thm}

{\bf Proof}\ \ Define: $A^*=\left(\begin{array}{cc}
0 & A\\
-A^T & 0
\end{array}
\right)=(a_{ij}^*)_{2n\times 2n}$ and $E^*=\{(i_l,n+j_l)|1\leq
l\leq k\}$. By Theorem 3.1, we have

$Pf(E^*(A^*))Pf(A^*)=Pf(E_p^*(A^*))Pf(\overline{E_p^*}(A^*))+$

$a^*_{i_p(n+j_p)}\sum\limits_{1\leq l\leq k, l\neq
p}a^*_{i_l(n+j_l)}
\left[f(p,l)Pf(E^*(A^*)_{\{i_p,n+j_l\}})Pf(A^*_{\{n+j_p,i_l\}})
-\right. $

$$\left.
g(p,l)Pf(E^*(A^*)_{\{i_p,i_l\}})Pf(A^*_{\{n+j_p,n+j_l\}})\right],\eqno{(4.1)}$$
where $f(p,l)=s([2n],i_p(n+j_l))s([2n],(n+j_p)i_l)$ and
$g(p,l)=s([2n],i_pi_l)s([2n],(n+j_p)(n+j_l))$. It is not difficult
to see that we have the following:
$$a^*_{i_p(n+j_p)}=a_{i_pj_p},\ a^*_{i_l(n+j_l)}=a_{i_lj_l},\ f(p,l)
=-(-1)^{i_p+j_p+i_l+j_l}.\eqno{(4.2)}$$ By Lemma 2.3 and the
definitions of $A^*$ and $E^*(A^*)$, we have
$$Pf(E^*(A^*))=(-1)^{\frac{1}{2}n(n-1)}\det(E[A]),\ Pf(A^*)=
(-1)^{\frac{1}{2}n(n-1)}\det(A), \eqno{(4.3)}$$
$$Pf(E_p^*(A^*))=(-1)^{\frac{1}{2}n(n-1)}\det(E_p[A]),\
Pf(\overline{E_p^*}(A^*))=(-1)^{\frac{1}{2}n(n-1)}\det(\overline{E_p}[A]).
\eqno{(4.4)}
$$
By $(i)$ in Lemma 4.5, we have
$$Pf(E^*(A^*)_{\{i_p,i_l\}})=0, \eqno{(4.5)}$$
and by $(ii)$ and $(iii)$ in Lemma 4.5, we have
$$Pf(E^*(A^*)_{\{i_p,n+j_l\}})=(-1)^{\frac{1}{2}(n-1)(n-2)}\det(E[A]_{i_pj_l}),
\eqno{(4.6)}$$$$
Pf(A^*_{\{n+j_p,i_l\}})=(-1)^{\frac{1}{2}(n-1)(n-2)}\det(A_{i_lj_p}).
\eqno{(4.7)}
$$
The theorem is immediate from $(4.1)-(4.7)$ and hence we have
completed the proof of the theorem. $\hfill\blacksquare$

In the proof of Theorem 4.2, $(4.1)$ is obtained from Theorem 3.1.
Obviously, A corresponding identity to $(4.1)$ can be obtained
from Corollary 3.2. Similarly, by this identity we can prove the
following:
\begin{thm}
Let $A=(a_{ij})_{n\times n}$ be a matrix of order $n$ and
$E=\{(i_l,j_l)|1\leq l\leq k\}$ a non empty subset of $[n]\times
[n]$, where $i_1\leq i_2\leq\ldots\leq i_k$. Then for a fixed
$p\in [k]$ we have
$$ \det(E[A])\det(A)$$
$$=\det(E_p[A]) \det(\overline{E_p}[A])-
\sum\limits_{1\leq l\leq k,l\neq p}(-1)^{i_l+j_l+i_p+j_p}
a_{i_pj_p}a_{i_lj_l}\det(E[A]_{i_lj_p})\det(A_{i_pj_l}), $$where
$E_p=E\backslash\{(i_p,j_p)\}$ and $\overline{E_p}=\{(i_p,j_p)\}$.
\end{thm}

The following result is immediate from Theorems 4.2 and 4.3.
\begin{cor}
Let $A=(a_{ij})_{n\times n}$ be a matrix of order $n$ and
$E=\{(i_l,j_l)|1\leq l\leq k\}$ a non empty subset of $[n]\times
[n]$, where $i_1\leq i_2\leq\ldots\leq i_k$. Then for a fixed
$p\in [k]$ we have
$$\sum\limits_{l=1}^k(-1)^{i_l+j_l}a_{i_lj_l}\{\det(A[E]_{i_pj_l})
\det(A_{i_lj_p})-\det(E[A]_{i_lj_p})\det(A_{i_pj_l})\}=0.$$
\end{cor}

\begin{ex}
Let $A=(a_{ij})_{3\times 3}, E=\{(1,1),(2,2),(3,3)\}$ and $p=2$.
Then, by Theorems 4.2 and 4.3, we have
$$ \left |\begin{array}{ccc}
a_{11} & a_{12} & a_{13}\\
a_{21} & a_{22} & a_{23}\\
a_{31} & a_{32} & a_{33}
\end{array}
\right| \left |\begin{array}{ccc}
0 & a_{12} & a_{13}\\
a_{21} & 0 & a_{23}\\
a_{31} & a_{32} & 0
\end{array}
\right|- \left |\begin{array}{ccc}
a_{11} & a_{12} & a_{13}\\
a_{21} & 0 & a_{23}\\
a_{31} & a_{32} & a_{33}
\end{array}
\right| \left |\begin{array}{ccc}
0 & a_{12} & a_{13}\\
a_{21} & a_{22} & a_{23}\\
a_{31} & a_{32} & 0
\end{array}
\right|
$$
$$=-a_{11}a_{22}\left |\begin{array}{cc}
a_{21} & a_{23}\\
a_{31} & a_{33}
\end{array}
\right| \left |\begin{array}{ccc}
a_{12} & a_{13}\\
a_{32} & 0
\end{array}
\right|-a_{22}a_{33}\left |\begin{array}{cc}
a_{11} & a_{13}\\
a_{21} & a_{23}
\end{array}
\right| \left |\begin{array}{ccc}
0 & a_{12}\\
a_{31} & a_{32}
\end{array}
\right|
$$
$$=-a_{11}a_{22}\left |\begin{array}{cc}
a_{12} & a_{13}\\
a_{32} & a_{33}
\end{array}
\right| \left |\begin{array}{ccc}
a_{21} & a_{23}\\
a_{31} & 0
\end{array}
\right|-a_{22}a_{33}\left |\begin{array}{cc}
a_{11} & a_{12}\\
a_{31} & a_{32}
\end{array}
\right| \left |\begin{array}{ccc}
0 & a_{13}\\
a_{21} & a_{23}
\end{array}
\right|.
$$
\end{ex}
\subsection{Graphical edge$-$condensation for enumerating
perfect matchings} \hspace*{\parindent}
Let $M(G)$ denote the sum of weights of perfect matchings of a
weighted graph $G$, where the weight of a perfect matching $M$ of
$G$ is defined as the product of weights of edges in $M$. It is
well known that computing $M(G)$ of a graph $G$ is an
$NP$-complete problem (see \cite{Jerr87}). Inspired by
$(2.2)$-Dodgson's Determinant$-$Evaluation Rule, Propp
\cite{Prop03} first found the method of graphical
vertex-condensation for enumerating perfect matchings of plane
bipartite graphs as follows:
\begin{prop}[Propp \cite{Prop03}]
Let $G=(U,V)$ be a plane bipartite graph in which $|U|=|V|$. Let
vertices $a, b, c$ and $d$ form a $4-$cycle face in $G$, $a, c\in
U$, and $b, d\in V$. Then
$$M(G)M(G-\{a,b,c,d\})=M(G-\{a,b\})M(G-\{c,d\})+M(G-\{a,d\})M(G-\{b,c\}).$$
\end{prop}

By a combinatorial method, Kuo \cite{Kuo04} generalized Propp's
result above as follows.
\begin{prop}[Kuo \cite{Kuo04}]
Let $G=(U,V)$ be a plane bipartite graph in which $|U|=|V|$. Let
vertices $a, b, c,$ and $d$ appear
in a cyclic order on a face of $G$.\\
{\bf (1)}\ \ If $a, c\in U$, and $b, d\in V$, then
$$M(G)M(G-\{a,b,c,d\})=M(G-\{a,b\})M(G-\{c,d\})+M(G-\{a,d\})M(G-\{b,c\}).$$
{\bf (2)}\ \ If $a, b\in U$, and $c, d\in V$, then
$$M(G-\{a,d\})M(G-\{b,c\})=M(G)M(G-\{a,b,c,d\})+M(G-\{a,c\})M(G-\{b,d\}).$$
\end{prop}

By Ciucu's Matching Factorization Theorem in \cite{Ciuc97}, Yan
and Zhang \cite{YZ05} obtained a more general result than Kuo's
for the method of graphical vertex-condensation for enumerating
perfect matchings of plane bipartite graphs. Furthermore, Yan et
al \cite{YYZ05} proved the following results:
\begin{prop}[Yan, Yeh and Zhang \cite{YYZ05}]
Let $G$ be a plane weighted graph with $2n$ vertices. Let vertices
$a_1, b_1, a_2, b_2, \ldots, a_k, b_k\ (2\leq k\leq n)$ appear in
a cyclic order on a face of $G$, and let $A=\{a_1, a_2, \cdots,
a_k\}$, $B=\{b_1, b_2, \cdots, b_k\}$. Then, for any $j=1, 2,
\cdots, k$, we have
$$\sum_{Y\subseteq B,\ |Y|\ \mbox {is odd}} M(G-a_j-Y)M(G-A\backslash \{a_j\}-\overline
Y)=\sum_{W\subseteq B,\ |W|\ \mbox{is even}}M(G-W)M(G-A-\overline
W), $$ where the first sum ranges over all odd subsets $Y$ of $B$
and the second sum ranges over all even subsets $W$ of $B$,
$\overline Y=B\backslash Y$ and $\overline W=B\backslash W$.
\end{prop}

The following result, which is a special case of the above
theorem, was first found by Kenyon and was sent to ``Domino Forum"
in an Email (for details, see \cite{YYZ05}).
\begin{cor}
Let $G$ be a plane graph with four vertices $a,b,c$ and $d$ (in
the cyclic order) adjacent to a single face. Then
$$
M(G)M(G-a-b-c-d)+M(G-a-c)M(G-b-d)
$$
$$
=M(G-a-b)M(G-c-d)+M(G-a-d)M(G-b-c).\eqno{(4.8)}$$
\end{cor}

By Ciucu's Matching Factorization Theorem in \cite{Ciuc97}, Yan et
al \cite{YYZ05} also obtained some results for the method of
graphical edge-condensation for enumerating perfect matchings of
plane graphs. In this subsection, by using the new Pfaffian
identity in Corollary 3.3 we will prove a quadratic relation,
which has a simpler form than the formula in \cite{YYZ05}, for the
method of graphical edge-condensation for computing perfect
matchings of plane graphs.

We first need to introduce the Pfaffian method for enumerating
perfect matchings \cite{Kast63,Kast67}. If $G^e$ is an orientation
of a simple graph $G$ and $C$ is a cycle of even length, we say
that $C$ is oddly oriented in $G^e$ if $C$ contains odd number of
edges that are directed in $G^e$ in the direction of each
orientation of $C$. We say that $G^e$ is a Pfaffian orientation of
$G$ if every nice cycle of even length of $G$ is oddly oriented in
$G^e$ (a cycle $C$ in $G$ is nice if $G-C$ has perfect matchings).
It is well known that if a graph $G$ contains no subdivision of
$K_{3,3}$ then $G$ has a Pfaffian orientation (see \cite{Litt75}).
McCuaig \cite{McCu}, McCuaig et al \cite{MRST97}, and Robertson et
al. \cite{RST99} found a polynomial-time algorithm to show whether
a bipartite graph has a Pfaffian orientation.
\begin{prop}[\cite{Kast67,Lova86}]
Let $G^e$ be a Pfaffian orientation of a graph $G$. Then
$$[M(G)]^2=\det(A(G^e),$$ where $A(G^e)$ is the skew adjacency matrix of $G^e$.
\end{prop}
\begin{rem}
Let $G^e$ be a Pfaffian orientation of a graph $G$ and $A(G^e)$
the skew adjacency matrix of $G^e$. By Cayley Theorem and
Proposition 4.5, we have
$$M(G)=\pm Pf(A(G^e)),$$
which implies that, for two arbitrary perfect matchings $\pi_1$
and $\pi_2$ of $G$, both $b_{\pi_1}$ and $b_{\pi_2}$ have the same
sign.
\end{rem}
\begin{prop}[Kasteleyn's theorem, \cite{Kast63,Kast67,Lova86}]
Every plane graph $G$ has an orientation $G^e$ such that every
boundary face-except possibly the unbounded face has an odd number
of edges oriented clockwise. Furthermore, such an orientation is a
Pfaffian orientation.
\end{prop}
Now we can prove the following result:
\begin{lem}
Let $G$ be a plane graph with four vertices $a,b,c$ and $d$ (in
the cyclic order) adjacent to the unbounded face. Let $G^e$ be an
arbitrary Pfaffian orientation satisfying the condition in
Proposition 4.6 and $A=A(G^e)$ the skew adjacency matrix of $G^e$.
Then all $Pf(A_{\{a,b,c,d\}})Pf(A),
Pf(A_{\{a,b\}})Pf(A_{\{c,d\}}),Pf(A_{\{a,c\}})Pf(A_{\{b,d\}})$ and
$Pf(A_{\{a,d\}})Pf(A_{\{b,c\}})$ have the same sign.
\end{lem}

{\bf Proof}\ \ By $(2.1)$ in Corollary 2.1, we have
$$Pf(A_{\{a,b,c,d\}})Pf(A)=
Pf(A_{\{a,b\}})Pf(A_{\{c,d\}})-Pf(A_{\{a,c\}})Pf(A_{\{b,d\}})
+Pf(A_{\{a,d\}})Pf(A_{\{b,c\}}). \eqno{(4.9)}$$ Obviously,
$A_{\{a,b,c,d\}}, A_{\{a,b\}}, A_{\{c,d\}}, A_{\{a,c\}},
A_{\{b,d\}}, A_{\{a,d\}}$ and $A_{\{b,c\}}$ are the skew adjacency
matrices of $G^e-a-b-c-d, G^e-a-b, G^e-c-d, G^e-a-c, G^e-b-d,
G^e-a-d$ and $G^e-b-c$, respectively. Note that all the
orientations $G^e-a-b-c-d, G^e-a-b, G^e-c-d, G^e-a-c, G^e-b-d,
G^e-a-d$ and $G^e-b-c$ of $G-a-b-c-d, G-a-b, G-c-d, G-a-c, G-b-d,
G-a-d$ and $G-b-c$ satisfy the condition in Proposition 4.6 and
hence are Pfaffian orientations. By Remark 4.3, we have
$$M(G)=\pm Pf(A),\ M(G-a-b-c-d)=\pm Pf(A_{\{a,b,c,d\}}).$$
Hence we have proved the following:
$$M(G)M(G-a-b-c-d)=\pm Pf(A)Pf(A_{\{a,b,c,d\}}). \eqno{(4.10)}$$
Similarly, we can prove the following:
$$M(G-a-b)M(G-c-d)=\pm Pf(A_{\{a,b\}})Pf(A_{\{c,d\}}); \eqno{(4.11)}$$
$$M(G-a-c)M(G-b-d)=\pm Pf(A_{\{a,c\}})Pf(A_{\{b,d\}}); \eqno{(4.12)}$$
$$M(G-a-d)M(G-b-c)=\pm Pf(A_{\{a,d\}})Pf(A_{\{b,c\}}). \eqno{(4.13)}$$
The lemma is immediate from $(4.8)-(4.13)$. $\hfill \blacksquare$

Now we can start to state the main result in this subsection as
follows.
\begin{thm}
Suppose $G$ is a plane weighted graph with even number of vertices
and the weight of every edge $e$ in $G$ is denoted by $\omega_e$.
Let $e_1=a_1b_1, e_2=a_2b_2, \ldots, e_k=a_kb_k$ ($k\geq 2$) be
$k$ independent edges in the boundary of a face $f$ of $G$, and
let vertices $a_1, b_1, a_2, b_2, \ldots, a_k, b_k$ appear in a
cyclic order on $f$ and let $X=\{e_i|\ i=1,2,\ldots,k\}$. Then,
for any $j=1,2,\ldots,k$,
$$M(G)M(G-X)=M(G-e_j)M(G-X\backslash\{e_{j}\})+$$
$$\omega_{e_j}\sum\limits_{1\leq i\leq k, i\neq
j}\omega_{e_i}[M(G-b_j-a_i)M(G-X-a_j-b_i)-M(G-b_j-b_i)M(G-X-a_j-a_i)].$$
\end{thm}

{\bf Proof}\ \ Note that $e_1=a_1b_1, e_2=a_2b_2, \ldots,
e_k=a_kb_k$ ($k\geq 2$) are $k$ independent edges in the boundary
of a face $f$ of $G$. It suffices to prove the following:
$$M(G)M(G-X)=M(G-e_1)M(G-X\backslash\{e_{1}\})+$$
$$\omega_{e_1}\sum\limits_{i=2}^k\omega_{e_i}[M(G-b_1-a_i)
M(G-X-a_1-b_i)-M(G-b_1-b_i)M(G-X-a_1-a_i)], \eqno{(4.14)}$$

Since $G$ is a plane graph, for an arbitrary face $F$ of $G$ there
exists a planar embedding of $G$ such that the face $F$ is the
unbounded one. Hence we may assume that vertices $a_1, b_1, a_2,
b_2, \ldots, a_k, b_k$ appear in a cyclic order on the unbounded
face of $G$. Let $T$ be a spanning trees containing $k$ edges
$e_i$'s and let $T^e$ be an orientation of $T$ such that the
direction of each edge $e_i$ is from $a_i$ to $b_i$ for
$i=1,2,\ldots,k$. Because each face of $G$ can be obtained from
$T$ by adding an edge, it is not difficult to see that there
exists an orientation $G^e$ of $G$ obtained from $T^e$ which
satisfies the condition in Proposition 4.6. Hence all $G^e, G^e-X,
G^e-e_j,G^e-X\backslash \{e_j\},G^e-a_i-b_j, G^e-X-a_j-b_i,
G^e-b_j-b_i$ and $G^e-X-a_i-a_j$ are Pfaffian orientations
satisfying the condition in Proposition 4.6, the skew adjacency
matrices of which are $A,
E(A),\overline{E_j}(A),E_j(A),A_{\{a_i,b_j\}},E(A)_{\{a_j,b_i\}},
A_{\{b_j,b_i\}}$ and $E(A)_{\{a_i,a_j\}}$, respectively, where
$E=\{(a_i,b_j)|1\leq i\leq k\}, E_j=E\backslash\{e_j\}$ and
$\overline{E_j}=E\backslash E_j$. Without loss of generality, we
may assume that $a_i=2i-1, b_i=2i$ for $i=1,2,\ldots,k$, that is,
$E=\{(1,2),(3,4),\ldots,(2k-1,2k)\}$. By Corollary 3.3, we have
$$Pf(E(A))Pf(A)=Pf(E_1(A))Pf(\overline{E_1}(A))+$$$$a_{12}
\sum\limits_{i=2}^ka_{2i-1,2i}\left[
Pf(E(A)_{\{1,2i\}})Pf(A_{\{2,2i-1\}})
-Pf(E(A)_{\{1,2i-1\}})Pf(A_{\{2,2i\}})\right].\eqno{(4.15)}$$ By a
similar method to that in Lemma 4.6, we can prove that
$$Pf(E(A))Pf(A)=\pm M(G-X)M(G); \eqno{(4.16)}$$
$$Pf(E_1(A))Pf(\overline{E_1}(A))=\pm M(G-X\backslash
\{e_1\})M(G-e_1); \eqno{(4.17)}$$
$$Pf(E(A)_{\{1,2i\}})Pf(A_{\{2,2i-1\}})=\pm
M(G-X-a_1-b_i)M(G-b_1-a_i);\eqno{(4.18)}$$
$$Pf(E(A)_{\{1,2i-1\}})Pf(A_{\{2,2i\}})=\pm
M(G-X-a_1-a_i)M(G-b_1-b_i). \eqno{(4.19)}$$ Since every perfect
matching of $G-X$ is also a perfect matching of $G$, by the
definition of the Pfaffian, both $Pf(A)$ and $Pf(E(A))$ have the
same sign. Hence by $(4.16)$ we have
$$Pf(E(A))Pf(A)=M(G-X)M(G). \eqno{(4.16')}$$
Similarly, we have
$$Pf(E_1(A))Pf(\overline{E_1}(A))=M(G-X\backslash
\{e_1\})M(G-e_1). \eqno{(4.17')}$$ Note that if $\pi'$ is a
perfect matching of $G-a_1-b_1-a_i-b_i$ ($i\neq 1$) then
$\pi=\pi'\cup \{(a_1,b_1),(a_i,b_i)\}$ is a perfect matching of
$G$. By the definition of the Pfaffian, it is not difficult to see
that both $b_{\pi}$ and $b_{\pi'}$ have the same sign, which
implies that both $Pf(A)$ and $Pf(A_{\{a_1,b_1,a_i,b_i\}})$ have
the same sign. Hence $Pf(A)Pf(A_{\{a_1,b_1,a_i,b_i\}})\geq 0$. By
Lemma 4.6, we have
$$Pf(A_{\{a_1,b_i\}})Pf(A_{\{b_1,a_i\}})\geq 0, Pf(A_{\{a_1,a_i\}})
Pf(A_{\{b_1,b_i\}})\geq 0. \eqno{(4.20)}$$ Since every perfect
matching of $G-X-a_1-b_i$ is also a perfect matching of
$G-a_1-b_i$, both $Pf(E(A)_{\{a_1,b_i\}})$ and
$Pf(A_{\{a_1,b_i\}})$ have the same sign. Similarly, both
$Pf(E(A)_{\{a_1,a_i\}})$ and $Pf(A_{\{a_1,a_i\}})$ have the same
sign. Hence by $(4.20)$ we have
$$Pf(E(A)_{\{a_1,b_i\}})Pf(A_{\{b_1,a_i\}})\geq 0, Pf(E(A)_{\{a_1,a_i\}})
Pf(A_{\{b_1,b_i\}})\geq 0. \eqno{(4.21)}$$ From $(4.18),(4.19)$
and $(4.21)$, we have
$$Pf(E(A)_{\{1,2i\}})Pf(A_{\{2,2i-1\}})=
M(G-X-a_1-b_i)M(G-b_1-a_i);\eqno{(4.18')}$$
$$Pf(E(A)_{\{1,2i-1\}})Pf(A_{\{2,2i\}})=
M(G-X-a_1-a_i)M(G-b_1-b_i). \eqno{(4.19')}$$

Note that $a_{12}=\omega_{e_1}$ and $a_{2i-1,2i}=\omega_{e_i}$. It
is not difficult to see that $(4.14)$ follows from $(4.15)$ and $
(4.16')-(4.19')$. Hence we have complete the proof of the theorem.
$\hfill \blacksquare$
\begin{rem}
The formula in Theorem 4.4 for the method of graphical
edge-condensation for enumerating perfect matchings of plane
graphs has a simpler form than that in Theorem 3.2 in \cite{YYZ05}
\end{rem}


\end{document}